\documentclass[letterpaper, 10 pt, conference]{ieeeconf}
\IEEEoverridecommandlockouts

\usepackage{amsmath}
\usepackage{graphicx}
\usepackage{amssymb}
\usepackage{color}
\usepackage{authblk}
\usepackage{cite}

\newtheorem{theorem}{Theorem}

\newtheorem{proposition}{Proposition}

\title{Model-Based Event-Triggered Control over Lossy Networks}

\author{Eloy Garcia and Panos J. Antsaklis
\thanks{E. Garcia is with the Control Science Center of Excellence, Air Force Research Laboratory, Wright-Patterson AFB, OH 45433. \ttfamily{eloy.garcia.2@us.af.mil}}
\thanks{Panos J. Antsaklis is with the Department of Electrical Engineering, University of Notre Dame, Notre Dame IN 46556.}
}

\begin{document}
\maketitle 

%\begin{keyword}                           
%   Cooperative control, Event-based control, Multi-agent systems, Consensus.                                 
%\end{keyword}                            

\begin{abstract}
The event-triggered control problem over lossy communication networks is addressed in this paper. Although packet dropouts have been considered in the implementation of event-triggered controllers, the assumption of protocols that employ acknowledgement messages persists. This paper provides an approach that relaxes such assumption. An event-based controller is implemented at the sensor node of the networked and uncertain system and it transmits feedback measurements to the controller node at asynchronous time instants. The transmitted packets of information are subject to dropouts by the lossy network. We show that the uncertain system can be asymptotically stabilized, that a positive minimum inter-event time exists, and that the proposed approach does not require acknowledgement messages.
\end{abstract}

\section{Introduction} \label{sec:one}

Event-triggered control of networked systems has received increased attention in recent years. The development of this control and scheduling approach has been motivated by the extensive use of digital communication networks with limited bandwidth. In networked systems, the communication channel is shared by different applications and in many instances only a reduced number of nodes are able to send information through the network within some specified time interval. The main goal in an event-triggered feedback control system is to reduce the number of instances where the sensor node needs to transmit feedback updates to the controller node.

In event-triggered broadcasting  \cite{AntaTabuada10,Astrom02,Donkers10,Garcia13,Tabuada07} a subsystem sends its local state to the network only when it is necessary, that is, only when a measure of the local subsystem state error is above a specified threshold. The problem of packet dropouts in event-triggered or similar frameworks has been addressed by different authors. 
In reference \cite{mamduhi2014event}, the authors considered packet losses for control of discrete-time systems. It was assumed that at every time step the scheduler is aware whether the transmission has been successful or not, which is equivalent to implementing acknowledgment (ACK) messages. Bommannavar and Basar \cite{bommannavar2008optimal} addressed a finite-time networked control problem similar to event-triggered control where there is a limit on the number of times that control signals can be transmitted to the plant. The authors also considered the case where transmissions are subject to packet dropouts.
In \cite{tallapragada2019event} the event-triggered control of linear systems with a similar constraint to packet dropouts was considered. The constraint is now given by channel blackouts, which are intervals of time during which no packet can be successfully transmitted.

Event-triggered control strategies have also been applied to stabilize multiple coupled subsystems as in \cite{Guinaldo11,Mazo11TAC,WangLemmon11,wang2010relaxing}. 
The paper \cite{guinaldo2014distributed}  considered delays and packet losses in the stabilization problem of coupled subsystems. Two event-triggered communication protocols were proposed in \cite{guinaldo2014distributed}. The first protocol preserves state consistency in the sense that all neighbors of a given node use the same version of the transmitted state by that node. In the presence of packet losses multiple retransmissions of the same measurement may be needed until it is guaranteed that all neighbors receive the update. At that point the transmitting node sends a permission message that allows neighbors to start using the new transmitted measurement. Due to disadvantages concerning retransmissions and associated overuse of the communication channel, a second protocol was presented  in \cite{guinaldo2014distributed} where the state consistency is relaxed and the neighbors of a given node are allowed to use different versions of the state of that node. However, acknowledgment (ACK) messages are required. In this case, retransmissions are used for neighbors whose ACK message is not received by the transmitting node. Meanwhile, the neighbors that successfully received the measurement update can use it to recalculate their control inputs regardless of any remaining neighbor nodes that have not received the update.

%The problem of event-triggered control over unreliable communication networks has been addressed by different authors. The paper \cite{guinaldo2014distributed} proposed two protocols for event-triggered control over unreliable communication networks. However, both of them relied on ACK signals to corroborate that transmitted packets are successfully received at their destination node.

Other approaches for event-triggered control over unreliable communication networks have been documented in \cite{yu2013model,dolk2015dynamic, Guinaldo12CDC,lehmann2012event}. However, the references above assume the implementation of reliable ACK messages. 
The assumption that ACK messages can be transmitted without packet dropouts makes the stability analysis easier by guaranteeing consistency of state error at different node locations.
In the present paper we consider the more general and also more practical case where ACK messages are not implemented since they may be subject to losses as well.

Seeking to relax the constrain imposed by ACK messages the authors of \cite{dolk2017event} provided an approach which does not employ the ACK scheme. This case presents an additional challenge since the transmitting device cannot distinguish between a successful and a failed transmission.
Reference \cite{dolk2017event} represents a rare instance where packet dropouts without ACK messages are considered in event-triggered control of networked systems over lossy networks. The work is motivated by established protocols that do not employ ACK messaging, such as, user datagram protocol (UDP). An additional motivation is that in lossy networks such as wireless communication networks, ACK messages can be dropped which will cause significant disruption to the event-based approaches that assume that ACK packets are never lost.
 
Similar to our work, reference  \cite{dolk2017event}  assumed the existence of a Maximum Allowable Number of Successive Dropouts (MANSD). One key aspect in  \cite{dolk2017event}   is that, in the implementation of the event-triggering mechanism, it is required to keep track of 2M+2 auxiliary variables which include past measurements that were transmitted to the controller node. The augmented state containing this variable could be significantly large depending on the value of the MANSD. In our work, we present a much simpler scheme where there is no need to augment the state with extra variables. A similar approach was considered in \cite{Garcia16iet}; however, events were generated  by implementing a combined state error and time-based threshold function. In other words, a maximum inter-event time was arbitrarily chosen and implemented in order to guarantee that a packet of information eventually arrives at the controller node. This is not the ideal case in event-triggered control since wasteful transmissions may be triggered when it is not necessary. Besides considering uncertain systems, this paper overcomes the limitations in \cite{Garcia16iet} by implementing a simpler threshold function which is based only on state errors.

The remainder of this paper is organized as follows. Section~\ref{sec:two} states the problem and describes the model-based event-triggered control approach. The analysis of state errors in the sensor and the controller nodes in the presence of packet dropouts is performed in  Section \ref{sec:Loss}. 
 Stability and Zeno behavior results are presented in Section \ref{sec:Results}. Examples are given in Section~\ref{sec:example} and Section~\ref{sec:conclusion} concludes the paper.

%%%%%%%%%%%%%%%%%%%%%%%%%%%%%%%%%%%%%%%%%%%%%%%%%%%%%%%%%%%%%%%%%%%%%%%%%%%%%%%
\section{Problem Statement} \label{sec:two}
%\textit{Notation}.
%The notations $1_N$ and $0_N$ represent column vectors of all ones and all zeros, respectively. $\mathbb{R}$ and $\mathbb{C}$ denote the set of real numbers and the set of complex numbers, respectively. For any $s\in\mathbb{C}$, $Re\left\{s\right\}$ represents the real part of $s$. $J_\nu^{\lambda_i}$ represents a Jordan block of size $\nu$ corresponding to eigenvalue $\lambda_i$ and $\otimes$ denotes the Kronecker product. The boldface $\textbf{e}^\lambda$ represents the exponential of the scalar $\lambda$ and $\textbf{e}^A$ represents the matrix exponential of matrix $A$.

\begin{figure}
	\begin{center}
		\includegraphics[width=8.4cm,trim=.4cm .4cm .4cm .4cm]{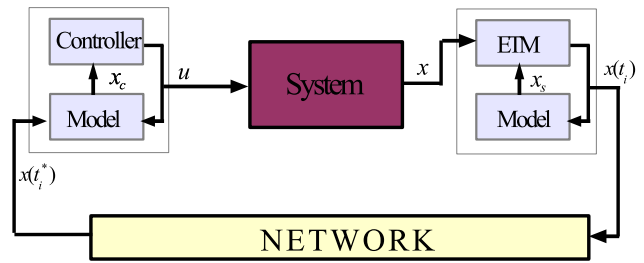}
	\caption{Model-based event-triggered control implementation}
	\label{fig:mbncs}
	\end{center}
\end{figure}

Consider the networked system  in Fig. \ref{fig:mbncs}
where a model-based networked control framework is implemented \cite{Garcia13scl,Garcia13MED}. 
We consider linear systems represented by
\begin{align}
  \left.
	\begin{array}{l l}
  \dot{x}(t) =Ax(t)+Bu(t)
\end{array}   \right.  \label{eq:plant}
\end{align}
where $x\in \mathbb{R}^n$ is the state of the system, $u\in \mathbb{R}^m$ is the control input, and the matrices $A$ and $B$ are uncertain matrices of appropriate dimensions. We assume that $A$ is unstable with distinct eigenvalues, with $n_u$ unstable eigenvalues and $n-n_u$ stable eigenvalues for $1 \leq n_u \leq n$. 

Since the system parameters are uncertain, then the available model matrices denoted by $\hat{A}$ and $\hat{B}$ are used for implementation of the controller and for evaluating the conditions to trigger events. Using $\hat{A}$ and $\hat{B}$, a model of the system is implemented in the sensor node and a second model of the system is implemented in the controller node. However, duet to packet losses, these two models do not generate the same estimate of the state. 
The dynamical model implemented at the sensor node is given by
\begin{align}
  \left.
	\begin{array}{l l}
  \dot{x}_s(t) =\hat{A}x_s(t)+\hat{B}u_s(t)  \\
  x_s(t_i)=x(t_i)
\end{array}   \right.  \label{eq:ModSensor}
\end{align}
where $u_s(t)=Kx_s(t)$ is the sensor model control input and $t_i$, for $i=0,1,2,...$, is the sequence of event time instants.

The dynamical model implemented at the controller node is given by
\begin{align}
  \left.
	\begin{array}{l l}
  \dot{x}_c(t) =\hat{A}x_c(t)+\hat{B}u(t)  \\
  x_c(t_{i^*})=x(t_{i^*})
\end{array}   \right.  \label{eq:ModCont}
\end{align}
where $u(t)=Kx_c(t)$ is the control input of the model implemented at the controller node. Note that $u(t)$, being the control input generated at the controller node, is also the control input of the real system \eqref{eq:plant}. The model based approach generalizes the zero-order-hold (ZOH) implementation by using the model parameters to estimate the state of the system between updates. The state error is reset when the controller node receives a state update. The model-based event-triggered (MB-ET) framework  helps to extend inter-event time intervals and, hence, reduce communication traffic, especially when the model parameters are a close representation of the real system parameters.

In the control architecture of Fig. \ref{fig:mbncs}  is important to note that, due to packet losses, the states of the models are not equal to each other all the time. Since ACK messages are not assumed in this paper, it is not possible to synchronize the state of the model in the sensor node with the state of the model in the controller node. The sensor model is updated more frequently than the controller model.  As it is shown in \eqref{eq:ModSensor},  the sensor model is updated every event time instant $t_i$ while the controller model is updated only at time instants $t_{i^*}$, when the transmitted update is not lost and it is successfully received at the controller node.

The sequence of time instants $t_{i^*}$ corresponds to the time instants when a packet of feedback information is successfully received at the controller node. The transmitted state is then used to update the state of the controller model \eqref{eq:ModCont}, that is, $x_c(t_{i^*})=x(t_{i^*})$.
 The sequence $t_{i^*}$ is a subset of the sequence $t_i$.
  %as shown in Fig. ???. For instance, in Fig.???, $i=0,1,2,3,...$ but $i^*=0,3,8,10...$. 
Then, we note that 
\begin{align}
	x_c(t)&=x_s(t), \ \ \  t\in[t_{i^*},t_{i^*+1} )   \\
	x_c(t)&\neq x_s(t), \ \ \  t\in[t_{i^*+1}, t_{i^*+M'}) 
\end{align}
where $1<M'\leq M$.

%The relationship between the two sequences of time instants is illustrated in Fig. ????.
The event time instants $t_i$ are generated by the event-triggering mechanism (ETM) located at the sensor node. Events are triggered according to  the condition
\begin{align}
	   t_{i+1} \!\!=\! \min \!\left\{t>t_i \big| \left\|e_s(t)\right\|> \beta \normalfont{\textbf{e}}^{-\alpha t} \right\} 
	    \label{eq:thre}
\end{align}
where  $\beta,\alpha>0$,  and 
\begin{align}
	  e_s(t)=x_s(t)-x(t).   \label{eq:errorS}
\end{align}
This means that the sequence $t_i$ represents the time instants at which events are generated at the sensor node because the state error $e_s(t)$ grows larger than the specified threshold. At these time instants an event is generated where the sensor model \eqref{eq:ModSensor} is updated using the state of the system $x_s(t_i)=x(t_i)$. With this update, the state error is reset to zero, that is, $e_s(t_i)=0$. Additionally, the state of the system is transmitted to the controller through the lossy communication channel. Some of these updates will be lost and will not arrive at the controller node. This is a particular problematic situation in event-triggered control; in contrast to periodic feedback, where an update is expected in periodic manner, the absence of updates in an event-triggered control implementation does not provide any useful information at the controller node. In other words, it is not known at the controller node whether an update was transmitted and lost or an update has not been transmitted at all.
In the following, we will analyze this problem. We will provide stability and Zeno exclusion results; but first, we analyze the response of the state of the model implemented at the controller node in the presence of packet dropouts.

%%%%%%%%%%%%%%%%%%%%%%%%%%%%%%%%%%%%%%%%%%%%%%%%%%%%%%%%%%%%%%%%%%%%%%%%%%%%%%%%%%%%%%%%%%%
\section{Packet Dropouts Analysis} \label{sec:Loss}
The assumption that ACK messages can be transmitted without packet dropouts makes the stability analysis easier by guaranteeing that the estimates $x_s(t)=x_c(t)$  for all $t>0$.
In this work we consider the more general and more practical case where ACK messages are not implemented since they may be subject to losses as well. In such a case the estimates $x_s(t)$ and $x_c(t)$ are not always equal to each other. 

Let us define the state error at the controller node
\begin{align}
	 e_c(t)=x_c(t)-x(t).   \label{eq:errorC}
\end{align}
Since $x_c(t)$ and $x_s(t)$ are different for some time intervals, then the error $e_c(t)$ and the error $e_s(t)$, defined in \eqref{eq:errorS}, are not the same. In general, the norm of the error $e_c(t)$ will grow larger than the norm of the error $e_s(t)$ since the model at the controller node is being updated less frequently than the model at the sensor node. The error $e_s(t)$, the error at the sensor node, can be continuously measured in order for the event scheduler to decide when to transmit the state. On the other hand, the error $e_c(t)$ represents the state error at the controller node and it cannot be measured since the state of the system, $x(t)$, is not available at the controller side of the networked system. How to manage and establish bounds on the error $e_c(t)$ in the presence of packet dropouts represents an important challenge. This problem is addressed in the remaining of this section.

Let us assume that there exist a Maximum Allowable Number of Successive Dropouts (MANSD) \cite{WangLemmon11,dolk2017event,ZhangYu08}, denoted as $M-1$, where $M>1$ is an integer. This means that if a measurement transmitted by the sensor node at time $t_{i^*}$ is successfully received at the controller node, then, at most $M-1$ consecutive dropouts are allowed and, in the worst case, the state measurement transmitted at time $t_{i^*+M}$ will be successfully received at the controller node.

In the following theorem we use the following notation to denote the current time in terms of the last successful update time instant. Let $t=t_{i^*}+\delta$, for $\delta\in[0,\delta_M)$. Also, let $0<\delta_j < \delta_{j+1}$ denote time instants such that $t_{i^*+j}=t_{i^*}+\delta_j$ and $t_{i^*+j+1}=t_{i^*}+\delta_{j+1}$, with $t_{i^*+j}<t_{i^*+j+1}$  for $j=1,...,M-1$. 

\begin{theorem}  \label{th:ContError}
Assume that the MANSD is $M-1$, for $M>1$. If the events at the sensor node are generated according to \eqref{eq:thre}
where $\beta,\alpha>0$,  
then, the error $e_c(t)$ at the controller node satisfies the following
\begin{align}
  \left.
	\begin{array}{l l}
   \left\|e_c(t)\right\|\!\!\! &\leq \Delta \beta \normalfont{\textbf{e}}^{-\alpha t}  
\end{array}   \right.  \label{eq:ContError}
\end{align}
for $t \geq 0$, where
\begin{align}
  \left.
	\begin{array}{l l}
  \Delta = 1 +   \sum_{k=1}^{M-1} \normalfont{\textbf{e}}^{\alpha \tilde{\delta}_k}|| \textbf{e}^{(\hat{A}+\hat{B}K)(\delta-\delta_k)} ||  
\end{array}   \label{eq:Normec} \right.
\end{align} 
and $\tilde{\delta}_k=\sum_{j=k}^{M-1} \bar{\delta}_j$. Also, $ \bar{\delta}_j = \frac{1}{\gamma+\alpha} \ln \big(\frac{\zeta_j+\beta\textbf{e}^{-\alpha t_{i^*+j}}}{\eta_j} \big)$ if $\kappa \geq \alpha$  and $ \bar{\delta}_j = \frac{1}{\gamma+\kappa} \ln \big(\frac{\zeta_j+\beta\textbf{e}^{-\alpha t_{i^*+j}}}{\eta_j} \big)$ if $\kappa < \alpha$, for positive constants $\gamma,\kappa,\eta_j>0$, and $\zeta_j\geq \eta_j$.
\end{theorem}
\textit{Proof.}
Let us consider the state error due to packet losses. Assume without loss of generality that the last update transmitted by the sensor node and successfully received at the controller node takes place at time instant $t_{i^*}$. 
%Hence, we have that $t_{k_i^j}=t_{k_i}$. 
The next event will be generated at time $t_{i^*+1}$; this event is generated because the condition \eqref{eq:thre} is satisfied. The current state of the system, $x(t_{i^*+1})$, will be transmitted from the sensor node to the controller node. Simultaneously, the state of the model at the sensor node is updated using the same state measurement that was just transmitted, that is, $x_s(t_{i^*+1})=x(t_{i^*+1})$. Thus, the error at the sensor node \eqref{eq:errorS} resets to zero, $e_s(t_{i^*+1})=0$. 

The error $e_s(t)$ did not trigger any event during the time interval $[t_{i^*},t_{i^*+1})$. By definition of the threshold function \eqref{eq:thre} and the reset $e_s(t_{i^*+1})=0$ at time $t_{i^*+1}$, we have that $\left\| e_s(t) \right\| = \left\| x_s(t)-x(t) \right\|  \leq \beta\normalfont{\textbf{e}}^{-\alpha t}$ for $t \in[t_{i^*},t_{i^*+1})$. Additionally, we have that 
\begin{align}
	\left\| e_s(t_{i^*+1}) \right\| = \left\| x_s(t^-_{i^*+1})-x(t_{i^*+1}) \right\|  = \beta\normalfont{\textbf{e}}^{-\alpha t_{i^*+1}}  \nonumber 
\end{align}
 where the notation $t^-_{i^*+1}$ denotes the time instant just before the update of the state of the model at the sensor node.  Note that we consider the state errors evaluated just before the event time instants $t_i$ (which can be denoted as $t_i^-$); however, to simplify notation, we will refer to the error $e(t)$ evaluated at time instants $t_i^-$ simply as $e(t_i)$.   

Let us now assume that the update at time $t_{i^*+1}$ is dropped, so $e_c(t_{i^*+1}) = \beta\textbf{e}^{-\alpha t_{i^*+1}}$. The next event will be generated at time $t_{i^*+2}$. The state of the system, $x(t_{i^*+2})$ will be transmitted to the controller node. Following similar steps we have that 
\begin{align}
\left\| e_s(t_{i^*+2}) \right\| = \left\| x_s(t^-_{i^*+2})-x(t_{i^*+2}) \right\|  = \beta\normalfont{\textbf{e}}^{-\alpha t_{i^*+2}} \nonumber
\end{align}
 holds. More generally, we have that
 \begin{align}
\left\| e_s(t_{i^*+j}) \right\| = \left\| x_s(t^-_{i^*+j})-x(t_{i^*+j}) \right\|  = \beta\normalfont{\textbf{e}}^{-\alpha t_{i^*+j}} \label{eq:SensErrj}
\end{align}
 for $j=1,...,M$.
Consider the worst case scenario where the number of successive dropouts after the last successful update at time $t_{i^*}$ is the MANSD, $M -1$. In such a case we have that the error $e_c(t)$ at time $t=t^-_{i^*+M}$, just before the update $x(t_{i^*+M})$ is successfully received at the controller node, satisfies the following

\begin{align}
  \left.
	\begin{array}{l l}
  \left\|e_c(t)\right\|\!\!\! &=||x_c(t)-x(t)||  \\
	  &=||x_c(t)-x(t)  \pm   x_s(t)  \\
	  &~~ \pm \textbf{e}^{(\hat{A}+\hat{B}K)(t-t_{i^*+M-1})}  x_s(t^-_{i^*+M-1})  \\  
	  &~~ \pm \textbf{e}^{(\hat{A}+\hat{B}K)(t-t_{i^*+M-2})}  x_s(t^-_{i^*+M-2})  \\  
				&~~~~ \vdots  \\
		  &~~ \pm \textbf{e}^{(\hat{A}+\hat{B}K)(t-t_{i^*+1})}  x_s(t^-_{i^*+1}) || 
\end{array}   \label{eq:Loss1} \right.
\end{align}   
for $t\in [t_{i^*+M-1},t_{i^*+M})$.  Given the sensor model dynamics and update law of the sensor model state defined in \eqref{eq:ModSensor}, we have that
\begin{align}
  \left.
	\begin{array}{l l}
       x_s(t^-_{i+1})  =  \textbf{e}^{(\hat{A}+\hat{B}K)(t_{i+1}-t_{i})}  x(t_{i})  
\end{array}    \right.  \nonumber  % \label{eq:SMResp}
\end{align}   
 for any event time instants $t_i$. Also note that, due to the controller model dynamics and update law of the controller model state defined in \eqref{eq:ModCont}, we have that
\begin{align}
  \left.
	\begin{array}{l l}
       x_c(t)  =  \textbf{e}^{(\hat{A}+\hat{B}K)(t-t_{i^*})}  x(t_{i^*})  
\end{array}   \right.  \nonumber
\end{align} 
for successful transmission time instants $t_{i^*}$. Thus, we can write \eqref{eq:Loss1} as follows
\begin{align}
  \left.
	\begin{array}{l l}
  \left\|e_c(t)\right\|  = \\
  ~~ || \textbf{e}^{(\hat{A}+\hat{B}K)(t-t_{i^*})}  x(t_{i^*}) -x(t) + x_s(t)  \\
	  ~~ + \textbf{e}^{(\hat{A}+\hat{B}K)(t-t_{i^*+M-1})}  [x_s(t^-_{i^*+M-1}) - x(t_{i^*+M-1})]  \\  
	  ~~ + \textbf{e}^{(\hat{A}+\hat{B}K)(t-t_{i^*+M-2})}  [x_s(t^-_{i^*+M-2}) - x(t_{i^*+M-2})]  \\  
				~~~~ \vdots  \\
		  ~~ + \textbf{e}^{(\hat{A}+\hat{B}K)(t-t_{i^*+1})}  [x_s(t^-_{i^*+1}) - x(t_{i^*+1})]  \\  
		 ~~ - \textbf{e}^{(\hat{A}+\hat{B}K)(t-t_{i^*})}  x(t_{i^*}) ||  \\
\end{array}   \label{eq:LossEv} \right.
\end{align}   
Cancelling the first and last term of the previous equation and writing the state differences in terms of the state error $e_s$, we have the following
\begin{align}
  \left.
	\begin{array}{l l}
  \left\|e_c(t)\right\| \!\!\! &=  ||  e_s(t)  \\
	&  ~~ + \textbf{e}^{(\hat{A}+\hat{B}K)(t-t_{i^*+M-1})}  e_s(t_{i^*+M-1})  \\  
	  &~~ + \textbf{e}^{(\hat{A}+\hat{B}K)(t-t_{i^*+M-2})}  e_s(t_{i^*+M-2})  \\  
		&		~~~~ \vdots  \\
		  &~~ + \textbf{e}^{(\hat{A}+\hat{B}K)(t-t_{i^*+1})}  e_s(t_{i^*+1}) ||.  \\
\end{array}   \label{eq:LossEv2} \right.
\end{align}   
We now note that each error term $||e_s(t_{i^*+j}) || =\beta\normalfont{\textbf{e}}^{-\alpha t_{i^*+j}} $ for $j=1,...,M-1$. Additionally, $||e_s(t) || \leq\beta\normalfont{\textbf{e}}^{-\alpha t}$; therefore  \eqref{eq:LossEv2}  satisfies
\begin{align}
  \left.
	\begin{array}{l l}
  \left\|e_c(t)\right\|  \!\!\! &\leq \beta\normalfont{\textbf{e}}^{-\alpha t} \\
  &~~+  \beta\normalfont{\textbf{e}}^{-\alpha t_{i^*+M-1}}|| \textbf{e}^{(\hat{A}+\hat{B}K)(t-t_{i^*+M-1})} ||  \\  
	&  ~~ + \beta\normalfont{\textbf{e}}^{-\alpha t_{i^*+M-2}}  || \textbf{e}^{(\hat{A}+\hat{B}K)(t-t_{i^*+M-2})} ||\\  
	&			~~~~ \vdots  \\
	&	  ~~ +\beta\normalfont{\textbf{e}}^{-\alpha t_{i^*+1}} || \textbf{e}^{(\hat{A}+\hat{B}K)(t-t_{i^*+1})} || .  \\
\end{array}   \label{eq:LossEv3} \right.
\end{align}   
%At this point a bound on the norm of the error $e_c$ has been obtained; however, this bound is given in terms of particular time instants
Using the notation described prior to this theorem, we have that $t=t_{i^*}+\delta$, where, in this case, $\delta\in[\delta_{M-1},\delta_M)$.  We can write the following expression
\begin{align}
  \left.
	\begin{array}{l l}
  \left\|e_c(t)\right\|  \!\!\! &\leq \beta\normalfont{\textbf{e}}^{-\alpha t}  \big( 1 +  \normalfont{\textbf{e}}^{\alpha (\delta-\delta_{M-1}})|| \textbf{e}^{(\hat{A}+\hat{B}K)(\delta-\delta_{M-1})} ||  \\  
	&  ~~ + \normalfont{\textbf{e}}^{\alpha(\delta- \delta_{M-2})}  || \textbf{e}^{(\hat{A}+\hat{B}K)(\delta-\delta_{M-2})} ||\\  
	&			~~~~ \vdots  \\
	&	  ~~ +\normalfont{\textbf{e}}^{\alpha (\delta-\delta_1)} || \textbf{e}^{(\hat{A}+\hat{B}K)(\delta-\delta_1)} || \big)  \\
\end{array}   \label{eq:LossEv4} \right.
\end{align}   
where $\delta-\delta_j>0$, for $j=1,...,M-1$. Let us now look at the response of the error $e_s(t)$ defined in \eqref{eq:errorS} in order to determine whether the differences $\delta-\delta_j$ are finite for $j=1,...,M-1$, that is, to determine if a successful feedback update eventually arrives at the controller node. Consider the following augmented system
\begin{align}
  \left.
	\begin{array}{l l}
   \begin{bmatrix}\dot{x} (t)\\ \dot{x}_c(t) \end{bmatrix} = \begin{bmatrix}A & BK \\ 0 & \hat{A}+\hat{B}K \end{bmatrix}   \begin{bmatrix} x (t) \\ x_c(t) \end{bmatrix}   
\end{array}   \label{eq:GammaSys} \right.
\end{align} 
with update law given by
\begin{align}
  \left.
	\begin{array}{l l}
     \begin{bmatrix} x (t_{i^*}) \\ x_c(t_{i^*}) \end{bmatrix}   =   \begin{bmatrix} x (t_{i^*}) \\ x(t_{i^*}) \end{bmatrix} .  
\end{array}   \label{eq:GammaUpd} \right.
\end{align} 
The update \eqref{eq:GammaUpd} simply denotes the fact that the state of the system is continuous. However, the state of the controller model is discontinuous and it is updated using $x(t_{i^*})$, that is, the model state is updated at time instants $t_{i^*}$ using the state of the system.
 Define 
\begin{align}
  \left.
	\begin{array}{l l}
  \Gamma = \begin{bmatrix}A & BK \\ \textbf{0}_n & \hat{A}+\hat{B}K \end{bmatrix}    \nonumber
\end{array}  \right.
\end{align} 
where the eigenvalues of $\Gamma$ are the eigenvalues of $A$ and the eigenvalues of $ \hat{A}+\hat{B}K$. Thus, $\Gamma$ is unstable since it contains the eigenvalues of $A$. Hence, $\Gamma$ has $n_u$ unstable eigenvalues and $2n-n_u$ stable eigenvalues for $1 \leq n_u \leq n$. Define $\mathcal{R}_s =span\{v_\iota\}$, for $\iota=n_u+1,...,2n$, where $v_\iota$ are the right eigenvectors of matrix $\Gamma$ associated with the stable eigenvalues $\lambda_\iota$, for $\iota=n_u+1,...,2n$.
%Also define $\mathcal{R}_u = \mathbb{R}^n - \mathcal{R}_s$.
Now, let us compute the response of the error $e_s$ as follows
\begin{align}
  \left.
	\begin{array}{l l}
  e_s(t) \!\!\!&= x_s(t)-x(t) \\
  &= \textbf{e}^{(\hat{A}+\hat{B}K)(t-t_{i^*+j})} x(t_{i^*+j})  - \Psi x(t_{i^*})
\end{array}   \label{eq:RespErrorS} \right.
\end{align}  
for $t\in [t_{i^*+j},t_{i^*+j+1})$ and $j=1,...,M-1$, where $e_s(t_{i^*+j})=0$ and 
\begin{align}
  \left.
	\begin{array}{l l}
  \Psi = [I_n \ \ \textbf{0}_n] \textbf{e}^{\Gamma(t-t_{i^*})}\begin{bmatrix}I_n \\ I_n \end{bmatrix} 
\end{array}   \label{eq:Psi} \right.
\end{align} 
where $I_n$ is the identity matrix of size $n$ and $\textbf{0}_n$ is a square matrix of size $n$ containing zeros in each entry.  

\begin{proposition}  \label{prop:unstable}
 The response of the state $x(t)$ can be bounded from below by an exponential term as follows
\begin{align}
  \left.
	\begin{array}{l l}
    ||x(t)|| = \left\| \Psi x(t_{i^*}) \right\| \geq \eta \textbf{e}^{\gamma (t-t_{i^*})}  
\end{array}   \right.   \nonumber
\end{align}
for  $t\in [t_{i^*},t_{i^*+M})$ and for  some $\eta, \gamma>0$. 
\end{proposition}
\textit{Proof}. The proof is provided in the Appendix. \ $\square$

We now proceed with the proof of the theorem.
The response of the state can be further described in terms of initial conditions at time $t_{i^*+j}$ as follows
\begin{align}
  \left.
	\begin{array}{l l}
 x(t) = \Psi_j \begin{bmatrix}x(t_{i^*+j}) \\ x_c(t_{i^*+j}) \end{bmatrix}   \nonumber
\end{array}   \right.
\end{align}
for $t\in [t_{i^*+j},t_{i^*+j+1})$ and $j=1,...,M-1$,  where
\begin{align}
 \Psi_j =[I_n \ \ \textbf{0}_n] \textbf{e}^{\Gamma(t-t_{i^*+j})}  \nonumber
\end{align}
and
\begin{align}
 \begin{bmatrix}x(t_{i^*+j}) \\ x_c(t_{i^*+j}) \end{bmatrix}  =  \textbf{e}^{\Gamma(t_{i^*+j}-t_{i^*})}\begin{bmatrix} I_n \\ I_n \end{bmatrix}  x(t_{i^*}).    \nonumber
\end{align}
The response of the state $x(t)$ can be bounded from below by an exponential term as follows
\begin{align}
  \left.
	\begin{array}{l l}
    ||x(t)|| = \left\| \Psi_j \begin{bmatrix}x(t_{i^*+j}) \\ x_c(t_{i^*+j}) \end{bmatrix} \right\| \geq \eta_j \textbf{e}^{\gamma (t-t_{i^*+j})}  
\end{array}   \right.   \nonumber
\end{align}
for  $t\in [t_{i^*+j},t_{i^*+j+1})$ and for  some $\eta_j>0$. 

On the other hand, the response of the model state $x_s$ in  \eqref{eq:ModSensor} is exponentially stable for $t\in [t_{i^*+j},t_{i^*+j+1})$ and there exist $\kappa>0$ and $\zeta_j \geq \eta_j$  such that
\begin{align}
  \left.
	\begin{array}{l l}
 || x_s(t) || \!\!\!& = || \textbf{e}^{(\hat{A}+\hat{B}K)(t-t_{i^*+j})} x(t_{i^*+j})  || \\
 &\leq  \zeta_j \textbf{e}^{-\kappa (t-t_{i^*+j})} 
\end{array}   \right.  \nonumber
\end{align}
for $t\in [t_{i^*+j},t_{i^*+j+1})$.
 Then, we have that
\begin{align}
  \left.
	\begin{array}{l l}
 || e_s(t) || \!\!\!& \geq \left\| \Psi_j \begin{bmatrix}x(t_{i^*+j}) \\ x_c(t_{i^*+j}) \end{bmatrix} \right\|  \\
 & \quad - || \textbf{e}^{(\hat{A}+\hat{B}K)(t-t_{i^*+j})} x(t_{i^*+j})  || \\
 &\geq \eta_j  \textbf{e}^{\gamma (t-t_{i^*+j})} -  \zeta_j\textbf{e}^{-\kappa (t-t_{i^*+j})} 
\end{array}    \right.  \nonumber
\end{align}
with $||e_s(t_{i^*+j})||=0$. Hence, the norm of $e_s(t)$ grows as time increases while the threshold function $\beta \textbf{e}^{-\alpha t}$ in  \eqref{eq:thre} decreases with time. Therefore, there exists a finite time instant $t_{i^*+j+1}>t_{i^*+j}$ such that $||e_s(t)||$ grows from zero at time instant $t_{i^*+j}$, to $\beta \textbf{e}^{-\alpha (t_{i^*+j+1})}$ at time instant $ t_{i^*+j+1}$. Thus, an upper-bound on the difference $t_{i^*+j+1}-t_{i^*+j}$ can be obtained by solving the following equation
\begin{align}
  \left.
	\begin{array}{l l}
 & \eta_j \textbf{e}^{\gamma (t_{i^*+j+1}-t_{i^*+j})} - \zeta_j \textbf{e}^{-\kappa (t_{i^*+j+1}-t_{i^*+j})} \\
 & \qquad =\beta \textbf{e}^{-\alpha t_{i^*+j+1}}.    \nonumber
\end{array}   \right.
\end{align}
Multiplying and dividing the right-hand side of the previous equation by $\textbf{e}^{\alpha t_{i^*+j}}$ we obtain
\begin{align}
  \left.
	\begin{array}{l l}
 & \eta_j \textbf{e}^{\gamma (t_{i^*+j+1}-t_{i^*+j})} - \zeta_j \textbf{e}^{-\kappa (t_{i^*+j+1}-t_{i^*+j})} \\
 & \qquad =\beta \textbf{e}^{-\alpha (t_{i^*+j+1}-t_{i^*+j})} \textbf{e}^{-\alpha t_{i^*+j}}  \\
  \Rightarrow &\eta_j \textbf{e}^{\gamma (\delta_{j+1}-\delta_j)} - \zeta_j \textbf{e}^{-\kappa (\delta_{j+1}-\delta_j)}  \\
  &\qquad =\beta \textbf{e}^{-\alpha (\delta_{j+1}-\delta_j)} \textbf{e}^{-\alpha t_{i^*+j}} .  \nonumber
\end{array}   \right.
\end{align}
In order to obtain an explicit upper-bound  $\bar{\delta}_j \geq \delta_{j+1}-\delta_j$ we consider two possible cases.  In the case where $\kappa \geq \alpha$ we have that 
\begin{align}
  \left.
	\begin{array}{l l}
  \eta_j \textbf{e}^{\gamma(\delta_{j+1}-\delta_j)} - \zeta_j \textbf{e}^{-\kappa (\delta_{j+1}-\delta_j)} \\
  \qquad\geq  \eta_j \textbf{e}^{\gamma (\delta_{j+1}-\delta_j)} - \zeta_j \textbf{e}^{-\alpha (\delta_{j+1}-\delta_j)} .  \nonumber
\end{array}   \right.
\end{align}
Then, the explicit upper-bound  $\bar{\delta}_j $ is obtained by solving the following equation
\begin{align}
  \left.
	\begin{array}{l l}
 \eta_j \textbf{e}^{\gamma \bar{\delta}_j} - \zeta_j \textbf{e}^{-\alpha \bar{\delta}_j}  =  \beta \textbf{e}^{-\alpha \bar{\delta}_j}  \textbf{e}^{-\alpha t_{i^*+j}}.  \nonumber
\end{array}   \right.
\end{align}
In this case the solution is given by
\begin{align}
  \left.
	\begin{array}{l l}
      \bar{\delta}_j = \frac{1}{\gamma+\alpha} \ln \big(\frac{\zeta_j+ \beta\textbf{e}^{-\alpha t_{i^*+j}}}{\eta_j} \big)   \label{eq:deltaa}
\end{array}   \right.
\end{align}
where $\frac{\zeta_j+\beta\textbf{e}^{-\alpha t_{i^*+j}}}{\eta_j}>1$, since $\zeta_j>\eta_j$. On the other hand, consider the case where $\kappa < \alpha$, then, we have that
\begin{align}
  \left.
	\begin{array}{l l}
  \beta \textbf{e}^{-\alpha \bar{\delta}_j}  <   \beta \textbf{e}^{-\kappa \bar{\delta}_j}   \nonumber
\end{array}   \right.
\end{align}
and the explicit upper-bound $\bar{\delta}_j$ in this case is obtained by solving the following equation
\begin{align}
  \left.
	\begin{array}{l l}
 \eta_j \textbf{e}^{\gamma \bar{\delta}_j} - \zeta_j \textbf{e}^{-\kappa \bar{\delta}_j}  =  \beta \textbf{e}^{-\kappa \bar{\delta}_j} \textbf{e}^{-\alpha t_{i^*+j}}.  \nonumber
\end{array}   \right.
\end{align}
In this case the solution is given by
\begin{align}
  \left.
	\begin{array}{l l}
      \bar{\delta}_j = \frac{1}{\gamma+\kappa} \ln \big(\frac{\zeta_j+\beta \textbf{e}^{-\alpha t_{i^*+j}}}{\eta_j} \big).  \label{eq:deltak}
\end{array}   \right.
\end{align} 
Finally, the divergent terms in \eqref{eq:LossEv4}, which are $\textbf{e}^{\alpha (\delta-\delta_k)}$, for $k=1,...,M-1$, can be upper-bounded by $\textbf{e}^{\alpha (\delta-\delta_k)}\leq \textbf{e}^{\alpha \tilde{\delta}_k}$ where $\tilde{\delta}_k=\sum_{j=k}^{M-1} \bar{\delta}_j$ and $\bar{\delta}_j$ is given by \eqref{eq:deltaa} and  \eqref{eq:deltak}. Also note that at each time instant $t_{i^*}$ we have that $e_c(t_{i^*})=0$. Therefore, the norm of the state error at the controller node is bounded by \eqref{eq:ContError} where $\Delta$ is given by \eqref{eq:Normec}.  $\square$

%Additionally, let us consider $\rho'<\rho$ and if some packet is received at $t_{i^*+\rho'}$ then we have that 
%\begin{align}
%  \beta\normalfont{\textbf{e}}^{\alpha \delta_{\rho'}}   \sum_{j=1}^{\rho'} \normalfont{\textbf{e}}^{-\alpha \delta_j}|| \textbf{e}^{(\hat{A}+\hat{B}K)(\delta_{\rho'}-\delta_j)} ||   \\
%  <  \beta\normalfont{\textbf{e}}^{\alpha \delta_\rho}   \sum_{j=1}^\rho \normalfont{\textbf{e}}^{-\alpha \delta_j}|| \textbf{e}^{(\hat{A}+\hat{B}K)(\delta_{\rho}-\delta_j)} || 
% \nonumber
%\end{align}
%since $\rho'<\rho$. Then, in general, eq. ???? holds

%????Analyze case  $x(t_i^*)\in \mathcal{R}_s$ ?????
%
%In Theorem \ref{th:ContError} it was assumed that $x(t_i^*)\in \mathcal{R}_u$, that is, the state vector transmitted from the sensor node belongs to the subspace  $\mathcal{R}_u = \mathbb{R}^n - \mathcal{R}_s$. In the case where $x(t_i^*)\in \mathcal{R}_s$ we have that 

%%%%%%%%%%%%%%%%%%%%%%%%%%%%%%%%%%%%%%%%%%%%%%%%%%%%%%%%%%%%%%%%%%%%%%%%%%%%%%%%%%%%%%%%%%%
\section{Stability and Minimum Inter-event Time} \label{sec:Results}
In this section we show that the networked system with events based on events triggered by \eqref{eq:thre} is stable under packet losses and without need for ACK messages. Let $x_0=||x(0)||$ and $\tilde{A}=A-\hat{A}$ and $\tilde{B}=B-\hat{B}$ denote the modeling error matrices. Also, we assume that the controller gain $K$ stabilizes the closed loop matrix $(A+BK)$. Then, there exist  constants $c,\bar{\alpha}>0$ such that $|| \textbf{e}^{(A+BK)t}|| \leq c\textbf{e}^{-\bar{\alpha}t}$. It is also shown in this section the existence of a positive inter-event time, that is, Zeno behavior is excluded.

\begin{theorem}\label{th:consensus}
Assume that the MANSD is $M-1$, for $M>1$. Then, system \eqref{eq:plant} with control input $u(t)=Kx_c(t)$, where $x_c(t)$ is generated by \eqref{eq:ModCont} is asymptotically stable in the presence of packet losses if the event time instants, $t_i$, are generated according to condition \eqref{eq:thre} where $\beta>0$ and $0<\alpha<\bar{\alpha}$.
Furthermore, the system does not exhibit Zeno behavior and the transmission inter-event times $t_{i+1}-t_i$  are bounded below by the \textit{positive} time $\underline{\delta}$, that is,
\begin{align}
 0< \underline{\delta}< t_{i+1}-t_i  \label{eq:tk}.
\end{align}	
The lower-bound on the inter-event times, $\underline{\delta}$, is explicitly given by
\begin{align}
    \left.
	\begin{array}{l l}
	\underline{\delta} = \frac{1}{\hat{a}+\bar{\alpha}}  \ln\big(1+ \frac{\beta}{\textbf{F}}\big)
	\end{array}  \label{eq:tau}   \right.
\end{align}
where $\hat{a}=||\hat{A}+\hat{B}K||$,
\begin{align}
\textbf{F}=\frac{\bar{F}\normalfont{\textbf{e}}^{(\alpha-\bar{\alpha}) t_i}}{\hat{a}+\bar{\alpha}}+ \frac{F }{\hat{a}+\alpha},   \label{eq:Htil}
\end{align}
$\bar{F}=\tilde{a}c(x_0 -\frac{\beta\Delta||BK||}{\bar{\alpha}-\alpha})$,  $F=(1+\frac{\tilde{a}c}{\bar{\alpha}-\alpha})\beta\Delta ||BK||$, $\tilde{a}=||\tilde{A}+\tilde{B}K||$, and $\Delta$ is given by \eqref{eq:Normec}.

\end{theorem}

\textit{Proof}.
We can write \eqref{eq:plant} as follows
\begin{align}
  \left.
	\begin{array}{l l}
  \dot{x}(t)  &=Ax(t) + BKx_c(t) \\
  &=Ax(t) + BK(x(t)+e_c(t)) \\
  &=(A+BK)x(t) + BKe_c(t).
\end{array}   \right.  \label{eq:plantSta}
\end{align}
Then, the response of the system is given by
\begin{align}
  x(t)=\textbf{e}^{(A+BK)t}x(0) + \int_0^t \textbf{e}^{(A+BK)(t-s)}BK e_c(s)ds    \nonumber
\end{align}
for $t \geq 0 $. Using the main result of the previous section given by the expression in \eqref{eq:ContError}, we can write
\begin{align}
  \left.
	\begin{array}{l l}
  || x(t) ||  = \left\| \textbf{e}^{(A+BK)t}x(0) +  \int_0^t \textbf{e}^{(A+BK)(t-s)}BK e_c(s)ds \right\|  \\
  \leq c \textbf{e}^{-\bar{\alpha}t} x_0 +  \int_0^t c \textbf{e}^{ -\bar{\alpha}(t-s)}||BK|| \Delta \beta  \textbf{e}^{-\alpha s} ds  \\
  \leq c  x_0 \textbf{e}^{-\bar{\alpha}t} + \frac{\beta c \Delta ||BK||}{\bar{\alpha}-\alpha} (\textbf{e}^{-\alpha t} -\textbf{e}^{-\bar{\alpha}t}).
\end{array}   \right.  \label{eq:normXt}
\end{align}
Since $0<\alpha<\bar{\alpha}$ we conclude that 
\begin{align}
\lim_{t\rightarrow\infty}  \left\| x(t)\right\| =0  \nonumber
\end{align}
and the system is asymptotically stable. 

In order to determine the lower-bound on positive minimum inter-event time, $\underline{\delta}$, we write the dynamics of the error at the sensor node as follows
\begin{align}
  \left.
	\begin{array}{l l}
  \dot{e}_s(t)& = \dot{x}_s(t) - \dot{x} (t) \\
  &= (\hat{A}+\hat{B}K)x_s(t) - Ax(t)-BK x_c(t)
  \end{array}   \right.  \label{eq:miet_errdyn}
\end{align}
for $t\in[t_i,t_{i+1})$, where $e_s(t_i)=0$.
 Substituting \eqref{eq:errorC} into \eqref{eq:miet_errdyn} we obtain
\begin{align}
  \left.
	\begin{array}{l l}
  \dot{e}_s(t) = (\hat{A}+\hat{B}K)e_s(t) - (\tilde{A}+\tilde{B}K) x(t) -BKe_c(t).
  \end{array}   \right.  \nonumber  % \label{eq:miet_errdyn}
\end{align}
Thus, we can write the following
\begin{align}
  \left.
	\begin{array}{l l}
 || \dot{e}_s(t)|| \leq  \hat{a}||e_s(t)|| +\tilde{a} ||x(t)|| +||BK|| ||e_c(t)||.
  \end{array}   \right.  \label{eq:miet_errdyn2}
\end{align}
We now substitute \eqref{eq:ContError} and \eqref{eq:normXt} into \eqref{eq:miet_errdyn2} to obtain
\begin{align}
  \left.
	\begin{array}{l l}
 || \dot{e}_s(t)|| \leq  \hat{a}||e_s(t)|| + \bar{F}\textbf{e}^{-\bar{\alpha} t}  + F \textbf{e}^{-\alpha t}
  \end{array}   \right.  \nonumber  %\label{eq:miet_errdyn3}
\end{align}
for $t\in[t_i,t_{i+1})$, where $||e_s(t_i)||=0$. Therefore, we have that the response of the error at the sensor node satisfies the following
\begin{align}
  \left.
	\begin{array}{l l}
 || e_s(t)|| & \leq  \frac{\bar{F}}{\hat{a}+\bar{\alpha}} (\textbf{e}^{\hat{a}t - (\hat{a}+\bar{\alpha}) t_i} - \textbf{e}^{-\bar{\alpha} t} ) \\
&~~+ \frac{F}{\hat{a}+\alpha} (\textbf{e}^{\hat{a}t - (\hat{a}+\alpha) t_i} - \textbf{e}^{-\alpha t} ).
  \end{array}   \right.  \nonumber  %\label{eq:miet_errdyn3}
\end{align}
Let us now substitute  $t=t_{i}+\delta$ into the previous equation and we have that
\begin{align}
  \left.
	\begin{array}{l l}
 || e_s(t)|| & \leq  \frac{\bar{F}}{\hat{a}+\bar{\alpha}} (\textbf{e}^{\hat{a} \delta} - \textbf{e}^{-\bar{\alpha} \delta} ) \textbf{e}^{-\bar{\alpha}t_i} \\
&~~+ \frac{F}{\hat{a}+\alpha} (\textbf{e}^{\hat{a}\delta} - \textbf{e}^{-\alpha \delta}) \textbf{e}^{-\alpha t_i} \triangleq f(\delta).
  \end{array}   \right.  \nonumber  %\label{eq:miet_errdyn3}
\end{align}
Therefore, the time $\delta>0$ that the function $f(\delta)$ takes to grow from zero, at time $t_i$, in order to reach the threshold $\beta\textbf{e}^{-\alpha t}=\beta\textbf{e}^{-\alpha(t_i+\delta)}$ is less or equal than the time it takes the error at the sensor node $\left\|e_s(t)\right\|$ to grow from zero, at time $t_i$, to reach the same threshold and generate the following event at time $t_{i+1}$, that is, $0<\delta \leq t_{i+1}-t_i$. Then, the minimum inter-event time is the greatest $\delta>0$ that satisfies the following expression
\begin{align}
    f(\delta) \leq \beta \textbf{e}^{-\alpha (\delta+t_i)}   \label{eq:miet_errdyn4}
\end{align}
which can also be written as follows
\begin{align}
  \left.
	\begin{array}{l l}
      \frac{\bar{F}}{\hat{a}+\bar{\alpha}} (\textbf{e}^{\hat{a}\delta} - \textbf{e}^{-\bar{\alpha} \delta} ) \textbf{e}^{(\alpha-\bar{\alpha})t_i} 
+ \frac{F}{\hat{a}+\alpha} (\textbf{e}^{\hat{a}\delta } - \textbf{e}^{-\alpha \delta})  \leq \beta \textbf{e}^{-\alpha \delta}.
  \end{array}   \right.  \nonumber  %\label{eq:miet_errdyn3}
\end{align}
An explicit expression for $\delta$, which is denoted by $\underline{\delta}$, can be obtained by noticing that
\begin{align}
  \left.
	\begin{array}{l l}
    \frac{\bar{F}}{\hat{a}+\bar{\alpha}} (\textbf{e}^{\hat{a}\delta} - \textbf{e}^{-\bar{\alpha} \delta} ) \textbf{e}^{(\alpha-\bar{\alpha})t_i} 
+ \frac{F}{\hat{a}+\alpha} (\textbf{e}^{\hat{a}\delta } - \textbf{e}^{-\alpha \delta})  \\
 \leq \frac{\bar{F}}{\hat{a}+\bar{\alpha}} (\textbf{e}^{\hat{a}\delta} - \textbf{e}^{-\bar{\alpha} \delta} ) \textbf{e}^{(\alpha-\bar{\alpha})t_i} 
+ \frac{F}{\hat{a}+\alpha} (\textbf{e}^{\hat{a}\delta } - \textbf{e}^{-\bar{\alpha} \delta}) 
  \end{array}   \right.  \nonumber  %\label{eq:miet_errdyn3}
\end{align}
and
\begin{align}
     \beta \textbf{e}^{-\alpha \delta} \geq  \beta \textbf{e}^{-\bar{\alpha} \delta}   \nonumber 
\end{align}
for $\delta\geq 0$, since $0<\alpha<\bar{\alpha}$. Thus, the solution of the following equation
\begin{align}
  \left.
	\begin{array}{l l}
   \frac{\bar{F}}{\hat{a}+\bar{\alpha}} (\textbf{e}^{\hat{a}\delta} - \textbf{e}^{-\bar{\alpha} \delta} ) \textbf{e}^{(\alpha-\bar{\alpha})t_i} 
+ \frac{F}{\hat{a}+\alpha} (\textbf{e}^{\hat{a}\delta } - \textbf{e}^{-\bar{\alpha} \delta})  =\beta \textbf{e}^{-\bar{\alpha} \delta}
  \end{array}   \right.  \nonumber  %\label{eq:miet_errdyn3}
\end{align}
guarantees that  \eqref{eq:miet_errdyn4} holds, that is, $\underline{\delta}<t_{i+1}-t_i$. The solution is given by $\underline{\delta}=\frac{1}{\hat{a}+\bar{\alpha}}  \ln\big(1+ \frac{\beta}{\textbf{F}}\big)$, which is clearly greater than zero, that is, the transmission inter-event time intervals are lower-bounded by a strictly positive value and the system does not exhibit Zeno behavior. 
 \ $\square$

%%%%%%%%%%%%%%%%%%%%%%%%%%%%%%%%%%%%%%%%%%%%%%%%%%%%%%%%%%%%%%%%%%%%%%%%%%%%%%%%%%%%5

\begin{figure}
	\begin{center}
		\includegraphics[height=9.2cm,width=9.6cm,trim=1.1cm .4cm -.4cm .4cm]{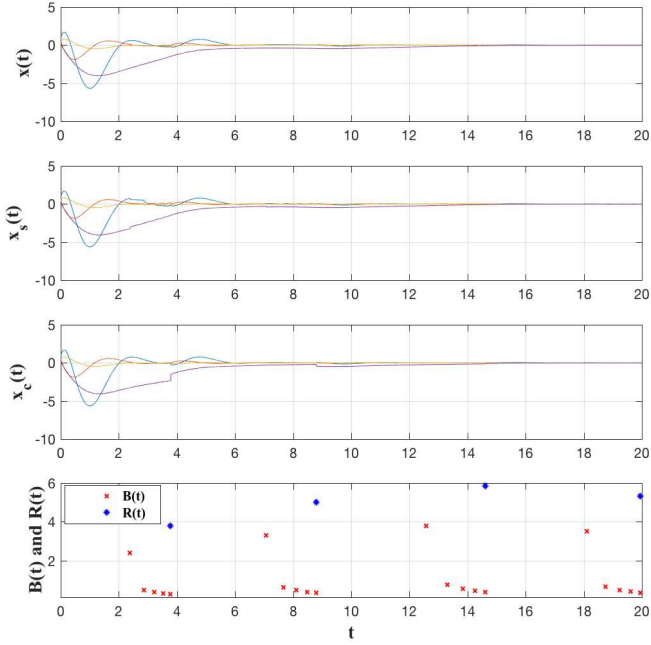}
	\caption{Response of the system and broadcasting and receiving time intervals for the MB-ET implementation}
	\label{fig:ex-mbet}
	\end{center}
\end{figure}

\section{Examples} \label{sec:example}
Consider the model of the vehicle for lane following control in \cite{Unyelioglu97} which is given by
\begin{align}
  \left.
	\begin{array}{l l}
    \begin{bmatrix}   \dot{V} \\ \dot{r} \\  \dot{\psi}  \\ \dot{Y}_g   \end{bmatrix}  \begin{bmatrix}   a_{11} &a_{12} & 0 & 0 \\    a_{21} &a_{22} & 0 & 0 \\  0 &1 & 0 & 0 \\ 1 & 0 & U & 0   \end{bmatrix}   \begin{bmatrix}   V \\ r \\  \psi  \\ Y_g   \end{bmatrix}  + \begin{bmatrix}   b_1 &0  \\ 0 &b_2 \\ 0 & 0  \\0 &0   \end{bmatrix} u  
\end{array}   \right.  \nonumber
\end{align}
where $V$ is the lateral velocity, $r$ is the yaw velocity, $\psi$ is the yaw position, and $Y_g$ is the $y$-axis coordinate of the vehicle's center of gravity. $U$ is the speed of the vehicle, considered in the opposite direction, and it is assumed to be constant. The steering control law is
\begin{align}
  \left.
	\begin{array}{l l}
 u =  \begin{bmatrix}   0 &0 & -K & K/d \\ 0 & 0 & -K & K/d   \end{bmatrix}  \begin{bmatrix}  \psi  \\ Y_g      \end{bmatrix}  
\end{array}   \right.  \nonumber
\end{align}
where $\psi$ and $Y_g$ can be used for event-triggered feedback and $d$ is the look-ahead distance. Consider the following nominal values $\hat{a}_{11}=-1.6579$, $\hat{a}_{12}= 10.4500$, $\hat{a}_{21}= 0.4886$,  $\hat{a}_{22}= -2.7180$, $\hat{b}_1 = -12.1053$, and ${b}_2 = 13.1429$. The real parameters contain perturbations around the nominal parameters $a=\hat{a}+\tilde{a}$ and $b=\hat{b}+\tilde{b}$, where $\tilde{a}\in[-0.1, 0.1]$ and $\tilde{b}\in[-0.05, 0.05]$. Let $d=40$, $KL=1$, and $U=-12$. In such a case the open-loop system is unstable and the aim is to stabilize it (converge to the lane and stay on it) using event-triggered feedback updates which may be subject to packet loss. Assume that $M=5$ and we choose the parameters $\beta=0.5$ and $\alpha=0.25$. 

The results of the simulation are shown in Fig. \ref{fig:ex-mbet}. The top plots show the state of the system $x(t)$ and the states of the models $x_s(t)$ and $x_c(t)$. It can be seen that, because of packet dropouts, $x_c(t)$, the state of the model at the controller node, is updated less frequently than $x_s(t)$, the state of the model at the sensor node.
The bottom plot shows the broadcasting/transmission time intervals and the receiving time intervals where we note that transmissions occur more frequently; however most of this transmitted packets are lost and they do not arrive at the controller node. The receiving time intervals measure the time interval between updates at the controller node, where a feedback measurement is successfully received.

\begin{figure}
	\begin{center}
		\includegraphics[width=8.0cm,trim=1.4cm .4cm 1.4cm .4cm]{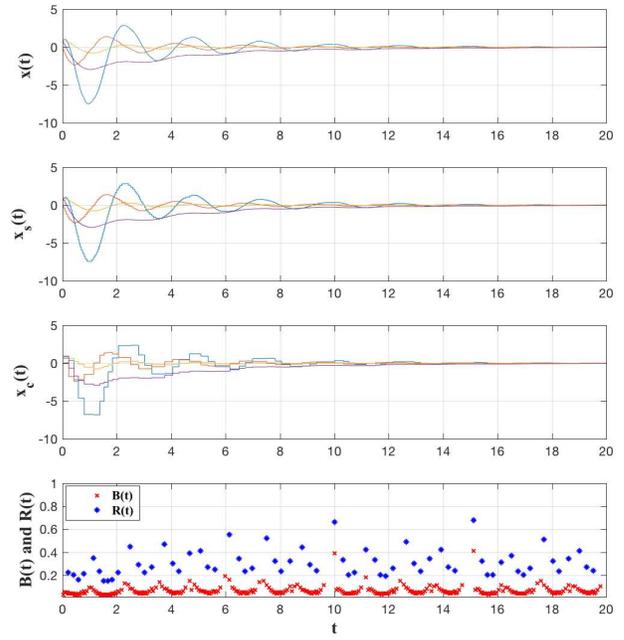}
	\caption{Response of the system and broadcasting and receiving time intervals for the ZOH-ET implementation}
	\label{fig:ex-zoh}
	\end{center}
\end{figure}

The approach discussed through this paper is not particular to the model-based event-triggered framework. The typical zero-order-hold (ZOH) event-triggered control approach also provides asymptotic stability of the system in the presence of packet losses and without ACK messages. 

\begin{theorem}  \label{th:ContErrorZOH}
Assume that the MANSD is $M-1$, for $M>1$; also assume that ZOH mechanisms are implemented at the sensor node and at the controller node. If the events at the sensor node are generated according to \eqref{eq:thre}
where $\beta,\alpha>0$,  
then, the error $e_c(t)$ at the controller  satisfies the following
\begin{align}
  \left.
	\begin{array}{l l}
   \left\|e_c(t)\right\|\!\!\! &\leq \Delta^{zoh} \beta \normalfont{\textbf{e}}^{-\alpha t}  
\end{array}   \right.  \label{eq:ContErrorZOH}
\end{align}
for $t \geq 0$, where
\begin{align}
  \left.
	\begin{array}{l l}
  \Delta^{zoh} =  \sum_{k=1}^{M} \normalfont{\textbf{e}}^{\alpha \tilde{\delta}^{zoh}_k}
  \end{array}   \label{eq:NormecZOH} \right.
\end{align} 
and $\tilde{\delta}^{zoh}_k=\sum_{j=k}^{M-1} \bar{\delta}^{zoh}_j$. Also, $ \bar{\delta}^{zoh}_j$ is the solution of the equation
\begin{align}
  \left.
	\begin{array}{l l}
   \Rightarrow &\eta^{zoh}_j \textbf{e}^{\gamma^{zoh} \bar{\delta}^{zoh}_j} - || x(t_{i^*+j})|| =\beta \textbf{e}^{-\alpha \bar{\delta}^{zoh}_j} \label{eq:zohnume}
\end{array}   \right.
\end{align} 
for $\gamma^{zoh}>0$, and $0<\eta^{zoh}_j \leq || x(t_{i^*+j})||$.
\end{theorem}
\textit{Proof.} The proof is given in Appendix B.

We will now compare the same example using the same event-triggered control parameters but instead of implementing a model-based approach, we use ZOH mechanisms. The results are  shown in Fig. \ref{fig:ex-zoh} where it can be seen that the system is stabilized; however, it takes longer for the system to settle at the steady-state response and, especially, the number of updates increases with respect to the MB-ET implementation presented in this paper.

\section{Conclusions} \label{sec:conclusion}
The stabilization of uncertain systems interconnected over an unreliable communication network subject to packet dropouts was addressed. The model-based event-triggered control framework was implemented together with the typical assumption of Maximum Allowable Number of Successive Dropouts. 
 The model-based event-triggered control framework extends inter-event time intervals, by using an estimate of the system state to generate the control input, compared to the zero-order-hold event-triggered control approach; this was illustrated in the examples. A great advantage is that this framework does not require the use of acknowledgement messages as it is common in event-triggered control approaches. This work provides a more general approach for control and stabilization over lossy communication networks.

\section*{Appendix A}   \label{sec:app}

\textit{Proof of Proposition \ref{prop:unstable}}. 
Even when the matrix $\Gamma$ is unstable it is not guaranteed that all trajectories of \eqref{eq:GammaSys} will be unstable. 
%Here, we look for the additional conditions so that all trajectories of \eqref{eq:GammaSys} with updates \eqref{eq:GammaUpd} grow exponentially.  Hence, it is guaranteed that the next event will be triggered in finite time when the repsonse of the system becomes unstable. 
Because $\Gamma$ has both stable and unstable eigenvalues, the unstable modes of $\Gamma$ can be suppressed if the update conditions of system \eqref{eq:GammaSys} are such that
\begin{align}
  \left.
	\begin{array}{l l}
      \begin{bmatrix} x (t_{i^*}) \\ x(t_{i^*}) \end{bmatrix} \in \mathcal{R}_s .
\end{array}   \right.   \label{eq:StableMode}
\end{align}
We will now show that \eqref{eq:StableMode} is not the case for system \eqref{eq:GammaSys} with updates given by \eqref{eq:GammaUpd}. This will guarantee that the next event will be triggered in finite time when the response of the system becomes unstable.

Let $\omega_\iota$ be the right eigenvector of matrix $A$ associated with its eigenvalue $\lambda_\iota$, for $\iota=1,...,n$. This means that $(\lambda_\iota I_n -A)\omega_\iota =0$. Also, let 
\begin{align}
  \left.
	\begin{array}{l l}
     v_\iota= \begin{bmatrix} \omega_\iota \\ 0 \end{bmatrix}
\end{array}   \right.   \label{eq:eigvectorA}
\end{align}
and compute 
\begin{align}
  \left.
	\begin{array}{l l}
     (\lambda_i I_{2n}-\Gamma)  v_\iota&=  \begin{bmatrix} \lambda_\iota I_n -A & -BK \\ \textbf{0}_n & \lambda_\iota I_n -(\hat{A}+\hat{B}K )    \end{bmatrix} \begin{bmatrix} \omega_\iota \\ 0 \end{bmatrix} \\
  & =  \begin{bmatrix} ( \lambda_\iota I_n -A )\omega_\iota \\ 0 \end{bmatrix} \\
  &= \begin{bmatrix}  0 \\ 0 \end{bmatrix}.
\end{array}   \right.   \nonumber 
\end{align}
Hence, the right eigenvectors of $\Gamma$ associated with the eigenvalues $\lambda_\iota$ of $A$ are given by \eqref{eq:eigvectorA}, for $\iota=1,...,n$.

Now, let $\chi_\iota$ be the right eigenvector of matrix $(\hat{A}+\hat{B}K ) $ associated with its eigenvalue $\lambda_\iota$, that is, $(\lambda_\iota I_n -(\hat{A}+\hat{B}K )  )\chi_\iota =0$. This second set of eigenvalues correspond to the eigenvalues $\lambda_\iota$ of $\Gamma$ for $\iota=n+1,...,2n$. Let us consider
\begin{align}
  \left.
	\begin{array}{l l}
     v_\iota= \begin{bmatrix}  \mu_\iota \\ \chi_\iota  \end{bmatrix}
\end{array}   \right.   \label{eq:eigvectorAhat}
\end{align}
and compute
\begin{align}
  \left.
	\begin{array}{l l}
     (\lambda_i I_{2n}-\Gamma)  v_\iota&=  \begin{bmatrix} \lambda_\iota I_n -A & -BK \\ \textbf{0}_n & \lambda_\iota I_n -(\hat{A}+\hat{B}K )    \end{bmatrix} \begin{bmatrix} \mu_\iota \\ \chi_\iota \end{bmatrix} \\
  & =  \begin{bmatrix} ( \lambda_\iota I_n -A )\mu_\iota -BK \chi_\iota  \\  (\lambda_\iota I_n -(\hat{A}+\hat{B}K )  )\chi_\iota \end{bmatrix}.
\end{array}   \right.   \nonumber 
\end{align}
Solving for $\mu_\iota$ in $( \lambda_\iota I_n -A )\mu_\iota -BK \chi_\iota =0$ we obtain 
\begin{align}
  \left.
	\begin{array}{l l}
     \mu_\iota &=  ( \lambda_\iota I_n -A )^{-1}BK \chi_\iota  \\
     & \neq \chi_\iota
\end{array}   \right.   \nonumber 
\end{align}
since  $( \lambda_\iota I_n -A )^{-1}BK \neq I_n$.

Hence, the right eigenvectors of $\Gamma$ associated with the eigenvalues $\lambda_\iota$ of $(\hat{A}+\hat{B}K )$ are given by \eqref{eq:eigvectorAhat}, for $\iota=n+1,...,2n$, where $  \mu_\iota  \neq \chi_\iota$. Thus, the set $\mathcal{R}_s$ is given by
\begin{align}
  \left.
	\begin{array}{l l}
\mathcal{R}_s = span \left\{ \begin{bmatrix}  \omega_\varsigma \\  0  \end{bmatrix}, \begin{bmatrix} ( \lambda_\iota I_n -A )^{-1}BK \chi_\iota   \\ \chi_\iota  \end{bmatrix} \right\}
\end{array}   \right.   \nonumber 
\end{align}
for $\omega_\varsigma$, for $\varsigma=n_u+1,...,n$, are the eigenvectors of $A$ associated with its stable eigenvalues and $\chi_\iota$, for $\iota=1,...,n$ are the eigenvectors of $(\hat{A}+\hat{B}K)$. Let us now write the following
\begin{align}
%  \left.
%	\begin{array}{l l}
  \sum_{\varsigma=n_u+1}^n a_\varsigma \begin{bmatrix}  \omega_\varsigma \\  0  \end{bmatrix} +  \sum_{\iota=1}^n b_\iota  \begin{bmatrix} ( \lambda_\iota I_n -A )^{-1}BK \chi_\iota   \\ \chi_\iota  \end{bmatrix}   \nonumber   \\
  = \begin{bmatrix}   \sum_{\varsigma=n_u+1}^n  a_\varsigma \omega_\varsigma  +  \sum_{\iota=1}^n b_\iota  ( \lambda_\iota I_n -A )^{-1}BK \chi_\iota   \\  \sum_{\iota=1}^n b_\iota   \chi_\iota \end{bmatrix}   \nonumber
%\end{array}   \right.   \nonumber 
\end{align}
for $a_\varsigma,b_\iota \in \mathbb{R}$, $\varsigma=n_u+1,...,n$, and $\iota=1,...,n$. Now, in order to determine whether \eqref{eq:StableMode} holds, we look at whether there exist $a_\varsigma$ such that 
\begin{align}
%  \left.
%	\begin{array}{l l}
  \sum_{\varsigma=n_u+1}^n  a_\varsigma \omega_\varsigma  +  \sum_{\iota=1}^n b_\iota  ( \lambda_\iota I_n -A )^{-1}BK \chi_\iota   = \sum_{\iota=1}^n b_\iota   \chi_\iota  \nonumber   \\
  \Rightarrow  \sum_{\varsigma=n_u+1}^n  a_\varsigma \omega_\varsigma  =  \sum_{\iota=1}^n b_\iota \big( I_n- ( \lambda_\iota I_n -A )^{-1}BK \big) \chi_\iota .   \label{eq:SysEq}  
%\end{array}   \right.   \nonumber 
\end{align}
Note that there are $n-n_u$ eigenvectors $\omega_\varsigma$, for $1\leq n_u \leq n$; hence, the previous equation results in an overdetermined system of equations and, in general, there are no solutions. 
%Define $\mathcal{R}_{As} =span\{\omega_\iota\}$, for $\iota=n_u+1,...,n$, where $\omega_\iota$ are the right eigenvectors of matrix $A$ associated with the stable eigenvalues $\lambda_\iota$, for $\iota=n_u+1,...,n$.
%Also define $\mathcal{R}_u = \mathbb{R}^n - \mathcal{R}_s$.

Then, we conclude that the updates \eqref{eq:GammaUpd} are such that 
\begin{align}
  \left.
	\begin{array}{l l}
      \begin{bmatrix} x (t_{i^*}) \\ x(t_{i^*}) \end{bmatrix} \notin \mathcal{R}_s 
\end{array}   \right.   \nonumber
\end{align}
and the response $x(t)$ grows exponentially. $\square$

%Finally, it is interesting to consider two particular cases. The first case is where all eigenvalues of $A$ are stable. In such a case we have that there are $n$ eigenvectors $\omega_\varsigma$ and there is a set of solutions for the system of equation \eqref{eq:SysEq}. The system of equations is no overdetermined any longer.  The second case is where there are no uncertainties, that is,  

\section*{Appendix B} \label{sec:appB}

\textit{Proof of Theorem \ref{th:ContErrorZOH}}. In the case where ZOH mechanisms are implemented in the sensor node and in the controller node, we have that
\begin{align}
	x_s(t)=x(t_i), \ \ \  t\in[t_i,t_{i+1} )   \nonumber 
\end{align}
and 
\begin{align}
	x_c(t)=x(t_{i^*}), \ \ \  t\in[t_{i^*},t_{i^*+M'} )   \nonumber 
\end{align}
for $1<M'\leq M$.
The analysis of the state error at the sensor node due to packet losses is similar to Theorem \ref{th:ContError}
and \eqref{eq:SensErrj} holds.

Consider the worst case scenario where the number of successive dropouts after the last successful update at time $t_{i^*}$ is the MANSD, $M -1$. In such a case we have that the error $e_c(t)$ at time $t=t^-_{i^*+M}$, just before the update $x(t_{i^*+M})$ is successfully received at the controller node, satisfies the following
\begin{align}
  \left.
	\begin{array}{l l}
  \left\|e_c(t)\right\|\!\!\! &=||x_c(t)-x(t)||  \\
	  &=||x(t_{i^*})-x(t)  \pm  x(t_{i^*+M-1})  \\  
	  &~~~~ \pm   x(t_{i^*+M-2}) \pm  ... \pm  x_s(t_{i^*+1}) ||  \\
	  &=||x(t_{i^*+M-1})-x(t)  \\
	  &~~~~ + [x(t_{i^*+M-2}) -x(t_{i^*+M-1})]  \\ 
	  &~~~~~~ \vdots \\
	   &~~~~ + [x(t_{i^*}) -x(t_{i^*-1})] || \\ 
\end{array}   \label{eq:Loss1ZOH} \right.
\end{align}   
for $t\in [t_{i^*+M-1},t_{i^*+M})$. Each state difference within brackets corresponds to the state error $e_s$, evaluated at different transmission time instants. Substituting  \eqref{eq:SensErrj}  into \eqref{eq:Loss1ZOH} we obtain 
\begin{align}
  \left.
	\begin{array}{l l}
  \left\|e_c(t)\right\|  \!\!\! &\leq \beta\normalfont{\textbf{e}}^{-\alpha t} +  \beta\normalfont{\textbf{e}}^{-\alpha t_{i^*+M-1}}  \\  
	&  ~~ + \beta\normalfont{\textbf{e}}^{-\alpha t_{i^*+M-2}}  +...  +\beta\normalfont{\textbf{e}}^{-\alpha t_{i^*+1}}.  
\end{array}   \label{eq:LossEv3ZOH} \right.
\end{align}   
We can now use the notation $t=t_{i^*}+\delta$, for $\delta\in[0,\delta_M)$ to obtain the following expression
\begin{align}
  \left.
	\begin{array}{l l}
 \left\|e_c(t)\right\|  \!\!\! &\leq \beta\normalfont{\textbf{e}}^{-\alpha t} [1 +  \normalfont{\textbf{e}}^{\alpha (\delta-\delta_{M-1})}  \\    
 &~~~~ +  \normalfont{\textbf{e}}^{\alpha (\delta-\delta_{M-2})}+...+  \normalfont{\textbf{e}}^{\alpha (\delta-\delta_1)}].
\end{array}   \label{eq:LossEv4ZOH} \right.
\end{align}   
Let us now look at the response of the error $e_s(t)$ defined in \eqref{eq:errorS} in order to determine whether the differences $\delta-\delta_j$ are finite for $j=1,...,M-1$, that is, to determine if a successful feedback update eventually arrives at the controller node. Consider the following augmented system
\begin{align}
  \left.
	\begin{array}{l l}
   \begin{bmatrix}\dot{x} (t)\\ \dot{x}_c(t) \end{bmatrix} = \begin{bmatrix}A & BK \\ \textbf{0}_n  & \textbf{0}_n  \end{bmatrix}   \begin{bmatrix} x (t) \\ x_c(t) \end{bmatrix}   
\end{array}   \label{eq:GammaSysZOH} \right.
\end{align} 
with update law given by
\begin{align}
  \left.
	\begin{array}{l l}
     \begin{bmatrix} x (t_{i^*}) \\ x_c(t_{i^*}) \end{bmatrix}   =   \begin{bmatrix} x (t_{i^*}) \\ x(t_{i^*}) \end{bmatrix} .  
\end{array}   \label{eq:GammaUpdZOH} \right.
\end{align} 
We now define 
\begin{align}
  \left.
	\begin{array}{l l}
  \Gamma^{zoh} = \begin{bmatrix}A & BK \\ \textbf{0}_n  & \textbf{0}_n  \end{bmatrix}    \nonumber
\end{array}  \right.
\end{align} 
where the first $n$ eigenvalues of $\Gamma$ are the eigenvalues of $A$ and the rest of its eigenvalues are equal to zero. Thus, $\Gamma$ is unstable since it contains the eigenvalues of $A$. Hence, $\Gamma$ has $n_u$ unstable eigenvalues and $2n-n_u$ either stable eigenvalues or eigenvalues at the origin, for $1 \leq n_u \leq n$. Define $\mathcal{R}_s =span\{v_\iota\}$, for $\iota=n_u+1,...,2n$, where $v_\iota$ are the right eigenvectors of matrix $\Gamma$ associated with the eigenvalues $\lambda_\iota$, for $\iota=n_u+1,...,2n$.
%Also define $\mathcal{R}_u = \mathbb{R}^n - \mathcal{R}_s$.
Now, let us compute the response of the error $e_s$ as follows
\begin{align}
  \left.
	\begin{array}{l l}
  e_s(t) \!\!\!&= x_s(t)-x(t) \\
  &=  x(t_{i^*+j})  - \Psi^{zoh} x(t_{i^*})
\end{array}   \label{eq:RespErrorSZOH} \right.
\end{align}  
for $t\in [t_{i^*+j},t_{i^*+j+1})$ and $j=1,...,M-1$, where $e_s(t_{i^*+j})=0$ and 
\begin{align}
  \left.
	\begin{array}{l l}
  \Psi^{zoh} = [I_n \ \ \textbf{0}_n ] \textbf{e}^{\Gamma^{zoh}(t-t_{i^*})}\begin{bmatrix}I_n \\ I_n \end{bmatrix} 
\end{array}   \label{eq:PsiZOH} \right.
\end{align} 
where $I_n$ is the identity matrix of size $n$. 

\begin{proposition}  \label{prop:unstablezoh}
 The response of the state $x(t)$ can be bounded from below by an exponential term as follows
\begin{align}
  \left.
	\begin{array}{l l}
    ||x(t)|| = \left\| \Psi^{zoh} x(t_{i^*}) \right\| \geq \eta^{zoh} \textbf{e}^{\gamma^{zoh} (t-t_{i^*})}  
\end{array}   \right.   \nonumber
\end{align}
for  $t\in [t_{i^*},t_{i^*+M})$ and for  some $\eta^{zoh}, \gamma^{zoh}>0$. 
\end{proposition}
\textit{Proof}. The proof is similar to the proof of Proposition \ref{prop:unstable} and it is ommitted. \ $\square$

Proceeding with the proof of the theorem we have that the response of the state can be further described in terms of initial conditions at time $t_{i^*+j}$ as follows
\begin{align}
  \left.
	\begin{array}{l l}
 x(t) = \Psi^{zoh}_j \begin{bmatrix}x(t_{i^*+j}) \\ x_c(t_{i^*+j}) \end{bmatrix}   \nonumber
\end{array}   \right.
\end{align}
for $t\in [t_{i^*+j},t_{i^*+j+1})$ and $j=1,...,M-1$,  where
\begin{align}
 \Psi^{zoh}_j =[I_n \ \ \textbf{0}_n ] \textbf{e}^{\Gamma^{zoh}(t-t_{i^*+j})}  \nonumber
\end{align}
and
\begin{align}
 \begin{bmatrix}x(t_{i^*+j}) \\ x_c(t_{i^*+j}) \end{bmatrix}  =  \textbf{e}^{\Gamma^{zoh}(t-t_{i^*+j})}\begin{bmatrix} I_n \\ I_n \end{bmatrix}  x(t_{i^*}).    \nonumber
\end{align}
The response of the state $x(t)$ can be bounded from below by an exponential term as follows
\begin{align}
  \left.
	\begin{array}{l l}
    ||x(t)|| = \left\| \Psi^{zoh}_j \begin{bmatrix}x(t_{i^*+j}) \\ x_c(t_{i^*+j}) \end{bmatrix} \right\| \geq \eta^{zoh}_j \textbf{e}^{\gamma^{zoh} (t-t_{i^*+j})}  
\end{array}   \right.   \nonumber
\end{align}
for  $t\in [t_{i^*+j},t_{i^*+j+1})$ and for  some $0<\eta^{zoh}_j \leq || x(t_{i^*+j})||$.

The norm of the error at the sensor node satisfies the following
\begin{align}
  \left.
	\begin{array}{l l}
 || e_s(t) || \!\!\!& \geq \left\| \Psi^{zoh}_j \begin{bmatrix}x(t_{i^*+j}) \\ x_c(t_{i^*+j}) \end{bmatrix} \right\| - || x(t_{i^*+j})|| \\
 &\geq \eta^{zoh}_j  \textbf{e}^{\gamma^{zoh} (t-t_{i^*+j})} - || x(t_{i^*+j})|| 
\end{array}    \right.  \nonumber
\end{align}
with $||e_s(t_{i^*+j})||=0$. Hence, the norm of $e_s(t)$ grows as time increases while the threshold function $\beta \textbf{e}^{-\alpha t}$ in  \eqref{eq:thre} decreases with time. Therefore, there exists a finite time instant $t_{i^*+j+1}>t_{i^*+j}$ such that $||e_s(t)||$ grows from zero at time $t_{i^*+j}$ to $\beta \textbf{e}^{-\alpha (t_{i^*+j+1})}$ at time $ t_{i^*+j+1}$. In addition, an upper-bound, $\bar{\delta}^{zoh}_j$, on the difference $t_{i^*+j+1}-t_{i^*+j}$ can be obtained by writing the following
\begin{align}
  \left.
	\begin{array}{l l}
% & \eta^{zoh}_j \textbf{e}^{\gamma^{zoh} (t_{i^*+j+1}-t_{i^*+j})} - || x(t_{i^*+j})||  \\
 %& \qquad =\beta \textbf{e}^{-\alpha (t_{i^*+j+1}-t_{i^*+j})} \\
  %\Rightarrow &\eta^{zoh}_j \textbf{e}^{\gamma^{zoh} (\delta_{j+1}-\delta_j)} - || x(t_{i^*+j})|| =\beta \textbf{e}^{-\alpha (\delta_{j+1}-\delta_j)} \\
   \Rightarrow &\eta^{zoh}_j \textbf{e}^{\gamma^{zoh} \bar{\delta}^{zoh}_j} - || x(t_{i^*+j})|| =\beta \textbf{e}^{-\alpha \bar{\delta}^{zoh}_j}.  \nonumber
\end{array}   \right.
\end{align}
The previous equation can be solved numerically for $\bar{\delta}^{zoh}_j$ where we note that $\bar{\delta}^{zoh}_j>0$ since $\eta^{zoh}_j \leq || x(t_{i^*+j})|| $.

%An upper-bound  $\bar{\delta}_j \geq \delta_{j+1}-\delta_j$ is obtained by solving the following equation
%\begin{align}
%  \left.
%	\begin{array}{l l}
% \eta^{zoh}_j \textbf{e}^{\gamma^{zoh} (\delta_{j+1}-\delta_j)} - || x(t_{i^*+j})|| =\beta \textbf{e}^{-\alpha (\delta_{j+1}-\delta_j)}.  \nonumber
%\end{array}   \right.
%\end{align}
%In this case the solution is given by
%\begin{align}
%  \left.
%	\begin{array}{l l}
%      \bar{\delta}_j = \frac{1}{\gamma+\alpha} \ln \big(\frac{\beta+\zeta_j}{\eta_j} \big)   \label{eq:deltaa}
%\end{array}   \right.
%\end{align}
%where $\frac{\beta+\zeta_j}{\eta_j}>1$. On the other hand, consider the case where $\kappa < \alpha$, then, we have that
%\begin{align}
%  \left.
%	\begin{array}{l l}
%  \beta \textbf{e}^{-\alpha \bar{\delta}_j}  <   \beta \textbf{e}^{-\kappa \bar{\delta}_j} .  \nonumber
%\end{array}   \right.
%\end{align}
%and the explicit upper-bound $\bar{\delta}_j$ in this case is obtained by solving the following equation
%\begin{align}
%  \left.
%	\begin{array}{l l}
% \eta_j \textbf{e}^{\gamma \bar{\delta}_j} - \zeta_j \textbf{e}^{-\kappa \bar{\delta}_j}  =  \beta \textbf{e}^{-\kappa \bar{\delta}_j}.  \nonumber
%\end{array}   \right.
%\end{align}
%In this case the solution is given by
%\begin{align}
%  \left.
%	\begin{array}{l l}
%      \bar{\delta}_j = \frac{1}{\gamma+\kappa} \ln \big(\frac{\beta+\zeta_j}{\eta_j} \big).  \label{eq:deltak}
%\end{array}   \right.
%\end{align} 
Finally, the divergent terms in \eqref{eq:LossEv4ZOH}, which are $\textbf{e}^{\alpha (\delta-\delta_k)}$, for $k=1,...,M-1$, can be upper-bounded by $\textbf{e}^{\alpha (\delta-\delta_k)}\leq \textbf{e}^{\alpha \tilde{\delta}^{zoh}_k}$ where $\tilde{\delta}^{zoh}_k=\sum_{j=k}^{M-1} \bar{\delta}^{zoh}_j$ and $\bar{\delta}^{zoh}_j>0$ is the solution of \eqref{eq:zohnume}. Also note that at each time instant $t_{i^*}$ we have that $e_c(t_{i^*})=0$. Therefore, the norm of the state error at the controller node is bounded by \eqref{eq:ContErrorZOH} where $\Delta^{zoh}$ is given by \eqref{eq:NormecZOH}.  $\square$

\bibliographystyle{IEEEtran}
\bibliography{ReferencesETLDc}

% Generated by IEEEtran.bst, version: 1.14 (2015/08/26)
\begin{thebibliography}{10}
\providecommand{\url}[1]{#1}
\csname url@samestyle\endcsname
\providecommand{\newblock}{\relax}
\providecommand{\bibinfo}[2]{#2}
\providecommand{\BIBentrySTDinterwordspacing}{\spaceskip=0pt\relax}
\providecommand{\BIBentryALTinterwordstretchfactor}{4}
\providecommand{\BIBentryALTinterwordspacing}{\spaceskip=\fontdimen2\font plus
\BIBentryALTinterwordstretchfactor\fontdimen3\font minus
  \fontdimen4\font\relax}
\providecommand{\BIBforeignlanguage}[2]{{%
\expandafter\ifx\csname l@#1\endcsname\relax
\typeout{** WARNING: IEEEtran.bst: No hyphenation pattern has been}%
\typeout{** loaded for the language `#1'. Using the pattern for}%
\typeout{** the default language instead.}%
\else
\language=\csname l@#1\endcsname
\fi
#2}}
\providecommand{\BIBdecl}{\relax}
\BIBdecl

\bibitem{AntaTabuada10}
A.~Anta and P.~Tabuada, ``To sample or not to sample: Self-triggered control
  for nonlinear systems,'' \emph{IEEE Transactions on Automatic Control},
  vol.~55, no.~9, pp. 2030--2042, 2010.

\bibitem{Astrom02}
K.~J. Astrom and B.~M. Bernhardson, ``Comparison of riemann and lebesgue
  sampling for first order stochastic systems,'' in \emph{41st IEEE Conference
  on Decision and Control}, 2002, pp. 2011--2016.

\bibitem{Donkers10}
M.~C.~F. Donkers and W.~P. M.~H. Heemels, ``Output-based event-triggered
  control with guaranteed linfty-gain and improved event-triggering,'' in
  \emph{49th IEEE Conference on Decision and Control}, 2010, pp. 3246--3251.

\bibitem{Garcia13}
E.~Garcia and P.~J. Antsaklis, ``Model-based event-triggered control for
  systems with quantization and time-varying network delays,'' \emph{IEEE
  Transactions on Automatic Control}, vol.~58, no.~2, pp. 422--434, 2013.

\bibitem{Tabuada07}
P.~Tabuada, ``Event-triggered real-time scheduling of stabilizing control
  tasks,'' \emph{IEEE Transactions on Automatic Control}, vol.~52, no.~9, pp.
  1680--1685, 2007.

\bibitem{mamduhi2014event}
M.~H. Mamduhi, D.~Toli{\'c}, A.~Molin, and S.~Hirche, ``Event-triggered
  scheduling for stochastic multi-loop networked control systems with packet
  dropouts,'' in \emph{53rd IEEE Conference on Decision and Control}.\hskip 1em
  plus 0.5em minus 0.4em\relax IEEE, 2014, pp. 2776--2782.

\bibitem{bommannavar2008optimal}
P.~Bommannavar and T.~Basar, ``Optimal control with limited control actions and
  lossy transmissions,'' in \emph{2008 47th IEEE conference on decision and
  control}.\hskip 1em plus 0.5em minus 0.4em\relax IEEE, 2008, pp. 2032--2037.

\bibitem{tallapragada2019event}
P.~Tallapragada, M.~Franceschetti, and J.~Cort{\'e}s, ``Event-triggered control
  under time-varying rates and channel blackouts,'' \emph{IFAC Journal of
  Systems and Control}, vol.~9, p. 100064, 2019.

\bibitem{Guinaldo11}
M.~Guinaldo, D.~V. Dimarogonas, K.~H. Johansson, J.~Sanchez, and S.~Dormido,
  ``Distributed event-based control for interconnected linear systems,'' in
  \emph{50th IEEE Conference on Decision and Control and European Control
  Conference (CDC-ECC)}, 2011, pp. 2553--2558.

\bibitem{Mazo11TAC}
M.~Mazo and P.~Tabuada, ``Decentralized event-triggered control over wireless
  sensor/actuator networks,'' \emph{IEEE Transactions on Automatic Control},
  vol.~56, no.~10, pp. 2456--2461, 2011.

\bibitem{WangLemmon11}
X.~Wang and M.~Lemmon, ``Event-triggering in distributed networked control
  systems,'' \emph{IEEE Transactions on Automatic Control}, vol.~56, no.~3, pp.
  586--601, 2011.

\bibitem{wang2010relaxing}
X.~Wang, Y.~Sun, and N.~Hovakimyan, ``Relaxing the consistency condition in
  distributed event-triggered networked control systems,'' in \emph{49th IEEE
  Conference on Decision and Control}, 2010, pp. 4727--4732.

\bibitem{guinaldo2014distributed}
M.~Guinaldo, D.~Lehmann, J.~S{\'a}nchez, S.~Dormido, and K.~H. Johansson,
  ``Distributed event-triggered control for non-reliable networks,''
  \emph{Journal of the Franklin Institute}, vol. 351, no.~12, pp. 5250--5273,
  2014.

\bibitem{yu2013model}
H.~Yu, E.~Garcia, and P.~J. Antsaklis, ``Model-based scheduling for networked
  control systems,'' in \emph{2013 American Control Conference}.\hskip 1em plus
  0.5em minus 0.4em\relax IEEE, 2013, pp. 2350--2355.

\bibitem{dolk2015dynamic}
V.~Dolk and W.~Heemels, ``Dynamic event-triggered control under packet losses:
  The case with acknowledgements,'' in \emph{2015 International Conference on
  Event-based Control, Communication, and Signal Processing (EBCCSP)}.\hskip
  1em plus 0.5em minus 0.4em\relax IEEE, 2015, pp. 1--7.

\bibitem{Guinaldo12CDC}
M.~Guinaldo, D.~Lehmann, J.~S{\'a}nchez, S.~Dormido, and K.~H. Johansson,
  ``Distributed event-triggered control with network delays and packet
  losses,'' in \emph{51st IEEE Conference on Decision and Control}, 2012, pp.
  1--6.

\bibitem{lehmann2012event}
D.~Lehmann and J.~Lunze, ``Event-based control with communication delays and
  packet losses,'' \emph{International Journal of Control}, vol.~85, no.~5, pp.
  563--577, 2012.

\bibitem{dolk2017event}
V.~Dolk and M.~Heemels, ``Event-triggered control systems under packet
  losses,'' \emph{Automatica}, vol.~80, pp. 143--155, 2017.

\bibitem{Garcia16iet}
{E. Garcia, Y. Cao, and D. W. Casbeer}, ``Decentralized event-triggered
  consensus of double integrator multi-agent systems with packet losses and
  communication delays,'' \emph{IET Control Theory \& Applications}, vol.~10,
  no.~15, pp. 1835--1843, 2016.

\bibitem{Garcia13scl}
E.~Garcia and P.~J. Antsaklis, ``Output feedback networked control with
  persistent disturbance attenuation,'' \emph{Systems \& Control Letters},
  vol.~62, no.~10, pp. 943--948, 2013.

\bibitem{Garcia13MED}
{E. Garcia and P. J. Antsaklis}, ``Model-based control of continuous-time
  systems with limited intermittent feedback,'' in \emph{21st Mediterranean
  Conference on Control and Automation}, 2013, pp. 452--457.

\bibitem{ZhangYu08}
W.~A. Zhang and L.~Yu, ``Modelling and control of networked control systems
  with both network induced delay and packet dropout,'' \emph{Automatica},
  vol.~44, no.~12, pp. 3206--3210, 2008.

\bibitem{Unyelioglu97}
K.~A. Unyelioglu, C.~Hatopoglu, and U.~Ozguner, ``Design and stability analysis
  of a lane following controller,'' \emph{IEEE Transactions on Control Systems
  Technology}, vol.~5, pp. 127--134, 1997.

\end{thebibliography}

\end{document}